\documentclass [12pt,a4paper] {article}
\usepackage[T1]{fontenc}
\usepackage{amssymb}
\usepackage[dvips]{color}
\usepackage{color}
\usepackage{amsmath,amsthm}
\usepackage{amsfonts}
\usepackage{amssymb}

\title{\color{mycolor4} Integral formulas for a Dirichlet series} 
\author{ Guy LAVILLE\\
 \\ Université de Caen\\
Laboratoire de Mathématiques Nicolas Oresme}
\date{  january 16,  2013 }
\definecolor{mycolor1}{rgb}{0.000,0.502,0.502}
\definecolor{mycolor2}{rgb}{0.502,0.000,0.502}
\definecolor{mycolor3}{rgb}{0.502,0.502,0.000}
\definecolor{mycolor4}{rgb}{0.502,0.000,0.000}
\newtheorem{theorem}{Theorem}
\def\cat#1{{\mathfrak{#1}}}

\begin{document}
\maketitle
\vskip 5. cm
\abstract{We present an integral representation formula for a Dirichlet series whose coefficients are the values of the Liouville's arithmetic function.  }
\newpage
\def\contentsname{Contents}
\tableofcontents 
\newpage
\section{Introduction }
    Let $ \sum_{n=1}^\infty \frac{a(n)}{n^s} $ be  a Dirichlet series such that : \\
- its analytic continuation is a meromorphic function with only one pole at $s=1$ \\
- there is a functional equation looking like :    
    \begin{displaymath}
    \sum_{n=1}^\infty \frac{a(n)}{n^s}=\varphi(s)\sum_{n=1}^{\infty}\frac{b(n)}{n^{1-s}} 
    \end{displaymath} \\
This gives a sequence $(b(n))$ allowing us to write a pseudo-cotangent or a pseudo-tangent function similar to a cotangent or a tangent function :    \\
\begin{displaymath}
\sum_{n=0}^{\infty} \frac{b(n)}{z^2+(2n+1)^2\pi^2}
\end{displaymath}
We prefer to choose a tangent function because the cotangent function has a singularity at the origin, hence some trouble to get a power series.   \\
It may be possible to deduce from the functional equation an integral formula for the starting Dirichlet series. Now we can hope to find a sequence $(c(n))$ allowing us to get an extension of this integral formula by a modification of the pseudo-tangent such as :     \\
\begin{displaymath}
\sum_{n=1}^\infty c(n)\frac{1}{\mathrm{e}^{z/n}+1}
\end{displaymath}
An easy example is the Riemann's $\zeta$ function itself, the three sequences are :  \\
\begin{displaymath}
(a(n))=(1,1,1,1,...) \qquad (b(n))=(1,1,1,1,...) \qquad (c(n))=(1,0,0,0,...)\\
\end{displaymath}
Another simple example is th Dirichlet series $ \sum_{n=1}^\infty \frac{(-1)^{n-1}}{n^s} $  :
\begin{displaymath}
(a(n))=(1,-1,1,-1,...) \qquad (b(n))=(1,0,1,0,...) \qquad (c(n))=(1,0,0,0,...)\\
\end{displaymath}
This is a general program. Here we take a particular case : a Dirichlet series equivalent to the Dirichlet series whose coefficients are the values of the Liouville arithmetic function. We obtain a representation integral formula. 
\section{Some functions associated with the Riemann's $ \zeta $ function }
\subsection{The functions $ \; \zeta, \; \zeta_a, \; \zeta_{imp} $ . }
The following Dirichlet functions are well known : 
\begin{align}
\zeta(s)&=\sum_{n=1}^\infty\frac{1}{n^s} \qquad \cat{R}(s)>1 \\
\zeta_a(s)&=\sum_{n=1}^\infty\frac{(-1)^{n-1}}{n^s} \qquad  \cat{R}(s)>1 \\
\zeta_{imp}(s)&=\sum_{m=0}^\infty\frac{1}{(2m+1)^s} \qquad \cat{R}(s)>1
\end{align}
The links with the $\zeta$  function are easy :
\begin{align}
\zeta(s)&=\frac{1}{1-2^{1-s}}\,\zeta_a(s)\\   
\zeta(s)&=\frac{1}{1-2^{-s}}\,\zeta_{imp}(s)   
\end{align}
Cf, for example [6] .
\subsection{The functions $ \; \zeta_{\lambda}, \; \zeta_\mu, \; \zeta_{\alpha}.$ }
Let $\lambda$ be the Liouville's arithmetic function :
$\lambda(1)=1$; for a prime $p$,  $ \lambda(p)=-1$; for all $a$ et $b$,   $\lambda(ab)=\lambda(a)\lambda(b)$. \\
  Let $\zeta_{\lambda} $ be the corresponding Dirichlet function : \begin{align}
\zeta_{\lambda}(s)&=\sum_{n=1}^\infty\frac{\lambda(n)}{n^s} \qquad \cat{R}(s)>1 \\
\zeta_{\lambda}(s)&=\frac{\zeta(2s)}{\zeta(s)}
\end{align}
$\zeta_{\lambda}$ is a meromorphic function on $\mathbb{C}$ . \\
Let $ \mu $ be the Möbius arithmetic fonction : \\
\begin{eqnarray}
\zeta_{\mu}(s) = \frac{1}{\zeta(s)} = \sum_{n=1}^\infty\frac{\mu(n)}{n^s} \qquad \cat{R}(s)>1
\end{eqnarray}
$\zeta_\lambda$ and $\zeta_\mu$ are also well known, cf [6] . \\
The singular points for $ \zeta_{\mu} $ are the zeros of $ \zeta $. But           $\zeta_{\lambda}$ has no singular points outside the domain  $  0 \leq \cat{R}(s) \leq 1 $. In that domain, its singular points are the zeros of the $\zeta$ function, except for $s=1$ . \\
Let $ \zeta_{\alpha} $ be :
\begin{eqnarray}
\zeta_{\alpha}(s)=\frac{\zeta_a(2s)}{\zeta_a(s)}
\end{eqnarray}
With (4) we get :
\begin{eqnarray}
\zeta_\lambda(s)=\frac{1-2^{1-s}}{1-2^{1-2s}}\zeta_{\alpha}(s) 
\end{eqnarray}
We do not need the arithmetic function $\alpha$ such that :
\begin{displaymath}
\zeta_\alpha(s)=\sum_{n=1}^\infty \frac{\alpha(n)}{n^s}
\end{displaymath}
\subsection{The function $\zeta_\beta$ . } 
Let $ \zeta_\beta $ be :
\begin{eqnarray}
\zeta_{\beta}(s) = \frac{\zeta_{imp}(2s-1)}{\zeta_{imp}(s)}
\end{eqnarray}
With (5) we have :
\begin{eqnarray}
\zeta_{\beta}(s)=\frac{(1-2^{1-2s})\zeta(2s-1)}{(1-2^{-s})\zeta(s)}
\end{eqnarray}
$\zeta_\beta $ is the generating function of an arithmetic function $ \beta $ :
\begin{equation}
\zeta_{\beta}(s) = \sum_{n=1}^{\infty} \frac{\beta(n)}{n^s} \quad \cat{R}(s)>1
\end{equation}
According to a theorem of Newman, cf [3], the following series is convergent and its value is :
\begin{equation}
\sum_{n=0}^\infty \frac{\mu(2n+1)}{2n+1} = 0
\end{equation}
Let $(2n+1)$ be an odd number. There is a unique decomposition in a factor without square and a square : 
\begin{eqnarray}
 \left\{
  \begin{array}{l l}
    2n+1=kh^2\\
    \beta(2n+1)=\mu(k)h \\
    \mid \beta(2n+1) \mid = h
  \end{array} \right.
\end{eqnarray}
We have the estimate :
\begin{eqnarray}
-1<\frac{\beta(2n+1)}{\sqrt{2n+1}}\leq 1
\end{eqnarray}
The equality is true if and only if $ 2n+1 $ is a square. \\
The Dirichlet series : 
\begin{eqnarray}
\zeta_{\beta}(s+\frac{1}{2}) = \sum_{n=0}^\infty \frac{\beta(n)}{n^{s+\frac{1}{2}}} = \sum_{n=1}^{\infty} \frac{\beta(n)}{\sqrt{n}} \frac{1}{n^s}
\end{eqnarray}
is, following (13) convergent for $ \cat{R}(s)> 1/2 $ . \\
The following inequality is useful :
\begin{equation}
\sum_{n=1}^\infty \frac{\mid \beta(n)\mid}{n^{3/2}} \leq \sum_{k=1}^\infty \frac{1}{k^{3/2}} \sum_{h=1}^\infty \frac{1}{h^2} 
\end{equation}
\subsection{The function $\zeta_\nu$ . } 
Let $ \zeta_\nu $ be :
\begin{eqnarray}
\zeta_\nu(s) = \frac{1}{\zeta_{imp}(s+1)} \frac{\zeta_{imp}(2s+2)}{\zeta_{imp}(s+3/2)} = \frac{\zeta_\beta(s+3/2)}{\zeta_{imp}(s+1)}
\end{eqnarray}
It is the generating function of an arithmetic function $\nu$ :
\begin{eqnarray}
\zeta_\nu(s) = \sum_{n=1}^\infty \frac{\nu(n)}{n^s} \quad \cat{R}(s)>0
\end{eqnarray}
Of course, $\nu$ is zero on even integers: $ \nu(2m)=0$ .
\begin{eqnarray}
\sum_{l\mid(2n+1)}l\nu(l) = \frac{\beta(2n+1)}{\sqrt{2n+1}}
\end{eqnarray}
the Möbius formula gives :
\begin{eqnarray}
(2n+1)\nu(2n+1) = \sum_{kl=2n+1} \mu(k) \frac{\beta(l)}{\sqrt{l}}
\end{eqnarray}
Let $d$ be the arithmetic function $d(n)=$ number of divisors of $n$. 
We have an estimate for all integers $m$ :
\begin{eqnarray}
\mid\nu(m)\mid \leq \frac{d(m)}{m}
\end{eqnarray}
\begin{theorem}
 The series whose terms are $\nu(n)$ converge and the value is $0$ .
\begin{eqnarray}
\sum_{n=1}^\infty\nu(n)=0.
\end{eqnarray}
\end{theorem}
Proof.
Use (22), (18) and (14). $\square$ \\
This result is very important for the present work .
\section{Integral formula for the $\zeta_a$ function } 
\subsection{The kernel $ \frac{1}{e^z+1} $ .}
The kernel $ \frac{1}{e^z+1} $ is better than $ \frac{1}{e^z-1} $  because we do not have a singularity at the origin. We have some classical expansions : 
\begin{align}
\frac{1}{e^t+1} &= \sum_{n=1}^\infty (-1)^{n-1}e^{-n t} \quad (t>0)   \notag \\
\frac{1}{e^z+1} &= \frac{1}{2}-2z\sum_{n=0}^\infty \frac{1}{z^2+(2n+1)^2\pi^2}   
\end{align}
Now, we take $ \mid z \mid < \pi $ for convergence .
\begin{align}
\frac{1}{e^z+1} &= \frac{1}{2}-2z\sum_{n=0}^\infty \frac{1}{(2n+1)^2\pi^2} \sum_{k=0}^\infty \frac{(-1)^k z^{2k}}{((2n+1)\pi)^{2k}} \\
\frac{1}{e^z+1} &= \frac{1}{2}-2\sum_{k=0}^\infty \frac{(-1)^k z^{2k+1}}
{\pi^{2k+2}} \zeta_{imp}(2k+2)
\end{align}
\subsection{Integral representation formula for $\zeta_a$ . }
It is well known that :
\begin{eqnarray}
\Gamma(s) \zeta_a(s) = \int_0^\infty \frac{1}{e^t+1} t^{s-1} dt \quad  \quad \cat{R}(s)>0
\end{eqnarray}
The continuation of this integral representation is possible by taking :  
\begin{displaymath}
\frac{1}{e^t+1}-\frac{1}{2}
\end{displaymath}
\section{Functional equations }
\subsection{Functional equation between $\zeta_a$ and $\zeta_{imp}$ .} 
For $\zeta$, the following functional equation of Riemann is well known cf [6] :
\begin{eqnarray}
\zeta(s) = 2^s \pi^{s-1} \sin\big(\frac{\pi}{2}s\big) \Gamma(1-s)\zeta(1-s)
\end{eqnarray}
thence a functional equation between $\zeta_a$ and $\zeta_{imp} $ :
\begin{eqnarray}
\zeta_a(s) = -2\pi^{s-1} \sin\big(\frac{\pi}{2}s\big) \Gamma(1-s) \zeta_{imp}(1-s)
\end{eqnarray}
\subsection{Functional equation between  $ \zeta_\alpha $ and $ \zeta_\beta $ .}
From (30) :
\begin{displaymath}
\frac{\zeta_a(2s)}{\zeta_a(s)} = \frac{-2\pi^{2s-1} \sin(\pi s) \Gamma(1-2s) \zeta_{imp}(1-2s)} {-2\pi^{s-1} \sin(\frac{\pi}{2}s) \Gamma(1-s) \zeta_{imp}(1-s)}
\end{displaymath}
Hence (cf [2] for the duplication formula) a first form for the functional equation between  $\zeta_\alpha$ and $\zeta_\beta$ :
\begin{eqnarray}
\zeta_\alpha(s) = 2^{1-2s} \; \pi^{s-1/2} \; \cos\big(\frac{\pi}{2}s\big)\; \Gamma\big(\frac{1}{2}-s\big) \; \zeta_\beta(1-s)
\end{eqnarray}
And a second form, but only for $ \cat{R}(s)<0 $ :
\begin{eqnarray}
\zeta_\alpha(s) = 2^{1-2s} \; \cos\big(\frac{\pi}{2}s\big)\; \Gamma\big(\frac{1}{2}-s\big) \; \sum_{m=0}^{\infty} \frac{\beta(2m+1)}{\sqrt{2m+1}} \frac{1}{(\pi (2m+1))^{1/2-s}}
\end{eqnarray}
\section{Integral formulas}
\subsection{Integral formula for $\zeta_\alpha$ in the domain $-3/2<\cat{R}(s)<-1/2$ . }
The functional equation (32) gives an integral for $\zeta_\alpha$ in the domain $ -3/2<\cat{R}(s)<0 $ : 
\begin{displaymath}
\zeta_\alpha(s) = 2^{1-2s} \cos\big(\frac{\pi}{2}s\big)\Gamma\big(\frac{1}{2}-s\big)\frac{2}{\pi} \cos\big(\frac{\pi}{2}s +\frac{\pi}{4} \big) \sum_{m=0}^{\infty} \frac{\beta(2m+1)}{\sqrt{2m+1}} \int_0^\infty \frac{x^{s+1/2}}{x^2+\pi^2 (2m+1)^2)} dx 
\end{displaymath}
Now, we want to permute the summation and the integral. To do this we have only to prove absolute integrability, but for $ -3/2<\cat{R}(s)<-1/2 $ . Let $ \sigma = \cat{R}(s) $, take the inequality (16)  : 
\begin{displaymath}
 \int_0^\infty 2x\sum_{n=1}^N  \left| \frac{\beta(n)}{\sqrt{n}(x^2+\pi^2 n^2)}  \right| \mid x^{s-1/2} \mid dx  
\leq \int_0^\infty \Big(\frac{1}{2}-\frac{1}{e^x+1} \Big) x^{\sigma -1/2} dx \quad -3/2<\cat{R}(s)<-1/2
\end{displaymath}
We get an integral formula for $\zeta_\alpha$ in $-3/2<\cat{R}(s)<-1/2$ :
\begin{eqnarray}
\zeta_\alpha(s) = \frac{2^{1-2s}}{\pi} \cos\big(\frac{\pi}{2}s\big) \cos\big(\frac{\pi}{2}s+\frac{\pi}{4}\big) \Gamma(1/2-s) \int_0^\infty 2x\sum_{m=0}^\infty \frac{\beta(2m+1)}{\sqrt{2m+1}(x^2+\pi^2(2m+1)^2} x^{s-1/2} dx  
\end{eqnarray}
Let :
\begin{displaymath}
\begin{split}
\varphi(s) &= \frac{2^{1-2s}}{\pi} \cos\big(\frac{\pi}{2}s\big) \cos\big(\frac{\pi}{2}s+\frac{\pi}{4}\big) \Gamma(1/2-s) \\
\mathcal{N}(x) &= 2x\sum_{m=0}^\infty \frac{\beta(2m+1)}{\sqrt{2m+1}\; (x^2+\pi^2(2m+1)^2)}
\end{split}
\end{displaymath}
Write (33) as :
\begin{eqnarray}
\zeta_\alpha(s) = \varphi(s) \int_0^\infty \mathcal{N}(x) x^{s-1/2} dx \quad -3/2<\cat{R}(s)<-1/2
\end{eqnarray}
Now, the aim is to prove this formula for a greater domain .
\subsection{The meromorphic function $\mathcal{N}$ .}
\begin{equation}
\mathcal{N}(z) = 2z\sum_{m=0}^\infty \frac{\beta(2m+1)}{\sqrt{2m+1}(z^2+\pi^2(2m+1)^2)}
\end{equation}
 is a meromorphic function in $\mathbb{C}$ .

All the poles are simple at $i\pi(2m+1)$ for $m \in \mathbb{Z} $. The residu is :  
\begin{displaymath}
\frac{\beta(2m+1)}{\sqrt{2m+1}}
\end{displaymath}
The expansion of $\mathcal{N}$ in a power series, in a neighborhood of zero is : 
\begin{align}
\mathcal{N}(z) &= 2 \sum_{k=0}^\infty \frac{(-1)^k z^{2k+1}}{\pi^{2k+2}} \zeta_\beta(2k+5/2) \qquad \mid z \mid <\pi
\end{align}
\subsection{Definition of the meromorphic function $ \mathcal{M} $ .}
Let $\mathcal{M}$ be the following meromorphic function in $\mathbb{C}$ .
\begin{eqnarray}
\mathcal{M}(z) =  \sum_{m=0}^\infty \nu(2m+1)\Big( \frac{1}{2} - \frac{1}{e^{z/(2m+1)}+1} \Big)
\end{eqnarray}
All the poles are simple at $ i(2l+1)\pi $ for $l \in \mathbb{Z} $.
The residu of $ \mathcal{M} $ is :
\begin{equation*}
\sum_{(2m+1)\mid (2l+1)} (2m+1)\nu(2m+1) = \frac{\beta(2l+1)}{\sqrt{2l+1}}
\end{equation*}
because of (21).  We obtain the same poles and the same residus. Of course, this does not give the equality betwwen $\mathcal{N}$ and $\mathcal{M}$.  \\
Theorem 1, (24) and the previous definition (37) of $ \mathcal{M} $ give us : 
\begin{equation}
\mathcal{M}(z) = -\sum_{m=0}^\infty \nu(2m+1) \frac{1}{e^{z/(2m+1)}+1}
\end{equation}
\subsection{Behavior of $ \mathcal{M} $ at infinity.} 
By Abel's summation by parts on (37), and theorem 1, (24), we get :
\begin{equation}
\lim_{x \rightarrow \infty} \quad \mathcal{M}(x) = 0
\end{equation}
\subsection{A bound for the derivative $\mathcal{M}' $ .} 
We can derive term by term the series of $\mathcal{M}(z) $ . From (38), we get :
\begin{equation}
\mathcal{M}'(z) = \sum_{m=0}^\infty \frac{\nu(2m+1)}{2m+1} \frac{e^{z/(2m+1)}}{(e^{z/(2m+1)}+1)^2}     
\end{equation}
Now, take this for $x$ real positive .
There exists a constant $C$ such that for all $ x \in [0,+\infty[ $ : 
\begin{equation}
\mid \mathcal{M}'(x) \mid \leq C 
\end{equation}
\subsection{A better bound for $ \mathcal{M} $ at infinity ?} 
Is it true that there exits a constant C such that : 
\begin{equation}
\mid \mathcal{M}(x) \mid \; \leq \frac{1}{x} C \quad \text{?}
\end{equation}
\subsection{Identity between $\mathcal{M}$ and $\mathcal{N}$ . }
The purpose is to prove that in a neighbohood of $0$, we have : 
$$ \mathcal{M}(z) = \mathcal{N}(z) $$
Starting from (27) we get for $\mid z \mid < \pi $ :
\begin{align*}
\sum_{n=0}^\infty \nu(2n+1) \Big( \frac{1}{2} - \frac{1}{e^{z/(2n+1)}+1} \Big) \quad &= \quad 2 \sum_{n=0}^\infty \nu(2n+1) \sum_{k=0}^\infty \frac{(-1)^k}{\pi^{2k+2}} \; \zeta_{imp}(2k+2) \; \frac{z^{2k+1}}{(2n+1)^{2k+1}}
\end{align*}
We can switch the sommations because we have absolute convergence.
Hence : 
\begin{align*}
\sum_{n=0}^\infty \nu(2n+1) \Big(\frac{1}{2} - \frac{1}{e^{z/(2n+1)}+1} \Big) =  2 \sum_{k=0}^\infty \frac{(-1)^k z^{2k+1}}{\pi^{2k+2}} \zeta_{imp}(2k+2) \zeta_\nu(2k+1) \\
\end{align*}
With (19), we get the power series of $\mathcal{M}$ at the origin : 
$$ \mathcal{M}(z) = 2 \sum_{k=0}^\infty \frac{(-1)^kz^{2k+1}} {\pi^{2k+1}} \zeta_\beta(2k+5/2) \qquad \mid z \mid < \pi $$
This is exactly (36) and we get that $ \mathcal{M}$ and $ \mathcal{N} $ are two expressions of the same function . 
\subsection{Integral representation formulas for $\zeta_\lambda$ in $-3/2<\cat{R}(s)<-1/2$ }
From (34) : 
$$ \zeta_\alpha(s) = \varphi(s) \int_0^\infty \mathcal{N}(x) x^{s-1/2} dx
\quad (-3/2<\cat{R}(s)<-1/2) $$
and with (42) we have : 
$$ \mid \mathcal{N}(x) \mid = \mid \mathcal{M}(x) \mid\leq C/x $$
Hence the convergence of the integral at $\infty$ for $\cat{R}(s)<1/2$ .
\begin{theorem}
In the domain $-3/2<\cat{R}(s)<-1/2$, we have the integral representation formula : 
\begin{equation}
\zeta_\lambda(s) = \frac{1-2^{1-s}}{1-2^{1-2s}} \; \varphi(s)\int_0^\infty \mathcal{N}(x)x^{s-1/2}dx
\end{equation}
\end{theorem}
More explicitely, with the values of $\varphi$, $\mathcal{N}$ and $\mathcal{M}$ :  
\begin{align}
\zeta_\lambda(s) &=  \frac{1-2^{1-s}}{(2^{2s-1}-1)\pi} \cos\big(\frac{\pi}{2}s\big) \cos\big(\frac{\pi}{2}s+\frac{\pi}{4}\big) \Gamma(1/2-s) \int_0^\infty \sum_{m=0}^\infty \frac{2x \; \beta(2m+1)} {\sqrt{2m+1}(x^2+\pi^2(2m+1)^2)} x^{s-1/2} dx \\
\zeta_\lambda(s) &= \frac{2^{1-s}-1}{(2^{2s-1}-1)\pi} \cos\big(\frac{\pi}{2}s\big) \cos\big(\frac{\pi}{2}s+\frac{\pi}{4}\big) \Gamma(1/2-s) \int_0^\infty \sum_{m=0}^\infty  \frac{\nu(2m+1)}{e^{x/(2m+1)}+1} x^{s-1/2} dx 
\end{align}
\section{References }
[1] A.E. Ingham. \textsl{The distribution of prime numbers}. Stechert-Hafner service agency, 1964.\newline
[2] W. Magnus, F. Oberhettinger, R.P. Soni. \textsl{Formulas and theorems for special functions of Mathematical Physics}. Springer-Verlag, 1966. \newline
[3] D.J. Newman. \textsl{Simple analytic proof of the prime number theorem}. Amer.Math.Monthly, 87(9) p.693-696, 1980 ou \textsl{Analytic Number theory}. Springer-Verlag, 1991. \newline
[4] B. Riemann. \textsl{Uber die Anzahl der Primzahlen unter einer gegebenen Grösse}. Monatsberichte der Berliner Akademie (1859), 671-680.\newline
[5] G. Tenenbaum. \textsl{Introduction à la théorie analytique et probabililiste des nombres}. Cours spécialisés. Société Mahtématique de France, 1995. \newline
[6] E.C. Titchmarsh. \textsl{The theory of the Riemann zeta-function}.  Oxford Science publication, second edition, 1986.
 \\
 \\
\textbf{Acknowledgements} \\  
I thank Michel Paugam for his many fruitful remarks and for his trust in my work, and also Claude Longuemare for what he taught to me through different subjects we studied.

\end{document}